\documentclass[leqno]{amsart}
\usepackage{amsmath}
\usepackage{amssymb}
\usepackage{amsthm}
\usepackage{enumerate}
\usepackage{graphicx}
\usepackage{verbatim}

\title{Surfaces and fronts 
with harmonic-mean curvature one   
       in hyperbolic three-space}
\author{Masatoshi Kokubu}
\date{\today}
\address{%
   Department of Natural Science,
   School of Engineering,
   Tokyo Denki University,
   2-2 Kanda-Nishiki-Cho,
   Chiyoda-Ku, Tokyo, 101-8457,
   Japan
}
\email{kokubu@cck.dendai.ac.jp}

\renewcommand\L{{\mathbb L}}
\renewcommand\H{{\mathbb H}}

\newcommand\la{\langle}
\newcommand\ra{\rangle}
\newcommand\iii{\sqrt{-1}}
\newcommand\fff{\, \mathrm{I}}
\newcommand\sff{\, \mathrm{I\!I}}
\newcommand\tff{\, \mathrm{I\!I\!I}}
\newcommand\pd{\partial}
\newcommand\R{{\mathbb R}}
\newcommand\C{{\mathbb C}}
\newcommand\D{{\mathbb D}}
\newcommand\lbar{\overline}
\newcommand{\SL}{\operatorname{SL}}
\newcommand{\SU}{\operatorname{SU}}

\newcommand{\PSL}{\operatorname{PSL}}
\newcommand{\Herm}{\operatorname{Herm}}
\newcommand{\trace}{\operatorname{tr}}
\newcommand\G{{\mathcal G}}
\newcommand\HH{{\mathcal H}}
\newcommand\W{{\mathcal W}}
\newcommand{\setdef}[3][;]{{\left\{#2\,#1\,#3\right\}}}
\newcommand\zb{\bar z}

\theoremstyle{plain}
\newtheorem{theorem}{Theorem}[section]
\newtheorem{proposition}[theorem]{Proposition}

\newtheorem{lemma}[theorem]{Lemma}
\newtheorem*{GMM}{G\'alvez-Mart\'\i{}nez-Mil\'an's formula}

\theoremstyle{definition}

\newtheorem*{remark}{Remark}
\newtheorem{example}{Example}

\numberwithin{equation}{section}

\begin{document}
\begin{abstract}
We investigate surfaces with constant harmonic-mean curvature one  
(HMC-1 surfaces) in hyperbolic three-space. 
We allow them to have certain kinds of singularities, 
and discuss some global properties.  
As well as flat surfaces and surfaces with constant mean curvature one
(CMC-1 surfaces),  
HMC-1 surfaces belong to a certain class of Weingarten surfaces. 
From the viewpoint of parallel surfaces, 
CMC-1 surfaces and HMC-1 surfaces are representative among 
this class.  

\end{abstract}

\maketitle

\section{Introduction}
In the differential geometry of surfaces in hyperbolic three-space $\H^3$, 
surfaces with constant mean curvature one 
(CMC-1 surfaces, for short) are one of the central subjects 
\cite{Br}, \cite{UY}, \cite{CHR}, etc.  
The theory of flat surfaces in $\H^3$ is also developing,  
thanks to the appearance of a representation formula 
due to G\'alvez, Mart\'\i{}nez and Mil\'an \cite{GMM1}. 
From the viewpoint of global theory for flat surfaces, 
one should generalize the 
category of surfaces to that of fronts. 
(Roughly speaking, a front is a surface with 
 certain kinds of singularities.) 
Any complete flat surface in $\H^3$ must be 
a horosphere or a hyperbolic cylinder, however, many complete flat 
fronts exist in $\H^3$ (see \cite{KUY}). 

On the other hand, G\'alvez, Mart\'\i{}nez and Mil\'an \cite{GMM2} 
also studied a wider class of surfaces in $\H^3$, including 
both CMC-1 surfaces and flat surfaces. 
It is a class of Weingarten surfaces satisfying $\alpha (H-1) = \beta K$ 
for some constants $\alpha$ and $\beta$. 
Here, $H$ denotes the mean curvature, and $K$ is the Gaussian curvature. 
Indeed, the following theorem is shown in \cite{GMM2}:

\begin{GMM}[\cite{GMM2}]
 Let $M$ be a non-compact, simply-connected surface and 
$f \colon M \to \H^3$ a Weingarten surface satisfying 
$\alpha(H-1)= \beta K$, where $\alpha$ and $\beta$ are real 
constants with $\alpha \ne 2 \beta$.
Then, there exist a meromorphic curve $\G \colon M \to \SL(2, \C)$ 
and a pair $(h, \theta)$ consisting of a meromorphic function $h$ 
and a holomorphic one-form $\theta$ on $M$, such that the immersion 
$f$ and its unit normal field $\nu$ can be recovered as  
$f = \G \HH \G^* \text{ and } \nu = \G \tilde \HH \G^*$, 
where
\begin{equation}\label{GMM:H}
\HH= \begin{bmatrix}
  \frac{1+\varepsilon^2 |h|^2}{1+ \varepsilon |h|^2} 
& -\varepsilon {\bar h} \\
-\varepsilon {h} & {1+\varepsilon |h|^2}
\end{bmatrix}
\text{ and } \quad 
\tilde \HH= \begin{bmatrix}
  \frac{1-\varepsilon^2 |h|^2}{1+ \varepsilon |h|^2} 
& \varepsilon {\bar h} \\
\varepsilon {h} & {-1-\varepsilon |h|^2}
\end{bmatrix}
\end{equation}
with $\varepsilon = \alpha /(\alpha - 2\beta)$ and 
$1+ \varepsilon |h|^2 >0$. 
Moreover, the curve $\G$ satisfies 
\begin{equation}\label{eq:g^-1dg}
 \G^{-1} d \G = \begin{bmatrix}
		 0 & \theta \\ dh & 0 
		\end{bmatrix}. 
\end{equation}
The following formulas hold: 
\begin{align}
 & \fff = (1-\varepsilon) \theta dh 
+ \left( \frac{(1-\varepsilon)^2 |dh|^2}{(1+\varepsilon |h|^2)^2}
 +(1+\varepsilon |h|^2)^2 |\theta|^2 \right)
+(1-\varepsilon) \bar \theta d \bar h, \label{GMM:fff}\\
 & \alpha \fff -2  \beta \sff = (\alpha - 2 \beta) \left(
(1+\varepsilon |h|^2)^2 |\theta|^2 
- \frac{(1-\varepsilon)^2 |dh|^2}{(1+\varepsilon |h|^2)^2}
\right), \label{GMM:sigma}
\end{align}
where $\fff$ and $\sff$ denote the first and second fundamental forms. 

Conversely, let $M$ be a Riemann surface, $\G \colon M \to \SL(2, \C)$ 
a meromorphic curve and $(h, \theta)$ a pair as above satisfying 
\eqref{eq:g^-1dg} and such that \eqref{GMM:sigma} is a positive 
definite metric. Then $f:= \G \HH \G^* \colon M \to \H^3$ 
{\rm ($\HH$ {\it as in} \eqref{GMM:H})}, is a Weingarten surface
 satisfying 
$\alpha (H-1) = \beta K$ 
with $\fff$ and $\alpha \fff -2  \beta \sff$ given 
by \eqref{GMM:fff} and \eqref{GMM:sigma}.
\end{GMM}
In the statement above, 
$\SL(2, \C)$ denotes the $2 \times 2$ complex special linear group, 
i.e., the complex Lie group consisting of 
$2 \times 2$ matrices with determinant 1, and we regard the hyperbolic 
$3$-space $\H^3$ as $\SL(2, \C)/{\SU(2)}$. 
(See Section \ref{sec:ov-GMM} for details.)

\bigskip
For the ratio ${[\alpha : \beta]}$ in $\R P^1$, let $\W_{[\alpha : \beta]}$ 
be the set of Weingarten surfaces satisfying 
$\alpha (H-1) = \beta K$, and set 
\begin{equation*}
 \W := \bigcup_{{[\alpha : \beta]} \in \R P^1} \W_{[\alpha : \beta]}. 
\end{equation*}
It is remarkable that $\W$ is closed under parallel transforms, that is, 
any parallel surface of any surface in $\W$ is always in $\W$. 
More precisely, 
dividing $\W$ into four subclasses
\begin{equation*}
\W^0:=\W_{[0:1]}, \quad \W^1:=\bigcup_{\lambda < 1/2} \W_{[1:\lambda]}, 
\quad \W^2:=\W_{[1:1/2]}, \quad
\W^3:=\bigcup_{\lambda > 1/2} \W_{[1:\lambda]},  
\end{equation*}   
we can prove that 
each $\W^j$ is closed under parallel transforms. 
(See Theorem \ref{th:para-1} and Theorem \ref{th:para-2}.)
Hence, we can roughly say that 
CMC-1 surface theory represents the theory of surfaces in 
$\W^1$. For instance, one can construct a Weingarten surface 
satisfying $H-1= \lambda K \ (\lambda < 1/2)$, 
though it may have singularities, 
by constructing any CMC-1 surface first 
and by parallelly transforming it appropriately. 
By the same reasoning, the theory of surfaces in $\W^3$ 
can be represented by one special type of surfaces. 
We will take $\W_{[1:1]}$ 
as that representative for $\W^3$, because surfaces in 
$\W_{[1:1]}$, i.e., Weingarten surfaces satisfying $H-1=K$, 
have another special geometric meaning: 
the sum of the reciprocals of the principal curvature is constantly 2. 
In other words, the harmonic mean of the principal curvature functions 
is constantly 1. 
We also call them {\it surfaces with constant harmonic-mean curvature one} 
(HMC-1 surfaces, for short).   

\bigskip

For the reason mentioned above, we will study HMC-1 surfaces 
in this paper. 
Although many works have been done on CMC-1 surfaces, 
HMC-1 surfaces have received less attention. 
For example, there is Epstein's work \cite{E}, however, it 
seems lesser-known.
(In classical Euclidean surface theory, 
the radii of principal curvature were considered as 
the fundamental entities. There seem to be 
some works about the mean radius of principal curvatures, or 
equivalently, about the harmonic mean of principal curvatures; e.g., 
Christoffel's theorem about rigidity of surfaces 
(cf. \cite[pp.299--302]{spivak}).)

In Section \ref{sec:bg}, we discuss background material for the sake of 
precisely understanding the contents mentioned in this introduction. 
Section \ref{sec:mrc-1surf} is devoted to deriving the formula 
due to G\'alvez, Mart\'\i{}nez and Mil\'an, but for the case of 
HMC-1 surfaces. 
Like the case of flat surfaces, it is more 
natural to consider HMC-1 fronts rather than HMC-1 surfaces. 
HMC-1 fronts are defined in Section \ref{sec:mrc-1front}.  
Some global properties are discussed and some examples are 
provided there.  

\medskip

The author would like to thank Professors Wayne Rossman, 
Masaaki Umehara and Kotaro Yamada for their valuable comment.  

\section{Background}\label{sec:bg}
\subsection{Basics}
Let $\L^4$ denote the Minkowski $4$-space with
the Lorentzian inner product $\la \, , \, \ra_L$ 
of signature $(-,+,+,+)$. 
Let $\mathcal F$ be the set of positively oriented and 
positively time-oriented frames 
$(e_0,e_1,e_2,e_3)$ in $\L^4$ satisfying 
\begin{equation}\label{eq:prod_e}
 \la e_{\alpha}, e_{\beta} \ra_L = 
\begin{cases}
 -1 & \text{ if } \alpha= \beta = 0, \\
 0 & \text{ if } \alpha \ne \beta, \\
 1 & \text{ if } \alpha= \beta > 0.  
\end{cases}
\end{equation}
The indices $\alpha$ and $\beta$ run over $0,1,2,3$, while 
the indices $i$, $j$ and $k$ run over $1,2,3$.  
  We shall use Einstein's convention, that is, the symbol $\sum$ 
may be omitted for sums over indices. 
  
Regarding $e_{\alpha} \colon (e_0,e_1,e_2,e_3) \ni \mathcal F \mapsto 
e_{\alpha} \in \L^{4} \ (\alpha = 0,1,2,3)$ as $\L^{4}$-valued functions, 
$d e_{\alpha} \ (\alpha = 0,1,2,3)$ are $\L^{4}$-valued one-forms 
on $\mathcal F$. The connection forms 
$\omega_{\alpha}^{\beta}$ are 
defined by $d e_{\alpha} =  e_{\beta} \otimes \omega_{\alpha}^{\beta}$. 
We write $\omega^i$ for $\omega_0^i$. 
Differentiating \eqref{eq:prod_e}, we have 
\begin{align}
&  \omega^{\alpha}_{\alpha}=0, \quad 
  -\omega_{i}^{0} + \omega_{0}^{i}=0, \quad 
  \omega_{i}^{j} + \omega^{i}_j =0, \label{eq:conn_form-1} \\
&  de_0 = e_{i} \otimes \omega^{i}, \quad \label{eq:de0ij}
 de_i = e_{0} \otimes \omega^{i} + e_{j} \otimes \omega_{i}^{j}.  
\end{align} 
Again, differentiating \eqref{eq:de0ij}, we have the structure equations:
\begin{equation}\label{eq:str-1}
  d \omega^i = -\omega^{i}_{j} \wedge \omega^j, \quad 
d \omega^i_j = -\omega^i_k \wedge \omega^k_j - \omega^i \wedge \omega^j. 
\end{equation}

The hyperbolic $3$-space $\H^3$ is the upper
half component of the two-sheeted hyperboloid in $\L^4$, 
i.e., 
\begin{equation*}
 \H^3 = \{x=(x_0,x_1,x_2,x_3)\in \mathbb L^4\,;\,
             \la x,x \ra_L=-1, \, x_0>0 \} 
\end{equation*}
with the metric induced by $\la \, , \, \ra_L$. $\H^3$ is 
a space form of constant negative curvature $-1$. As usual, 
we regard $e_0 \colon \mathcal F \to \H^3 \subset \L^4$ as 
the oriented orthonormal frame bundle of $\H^3$. 

Let $M$ be a connected, oriented surface, and 
$f\colon M \to \H^3$ an immersion. 
Let $\{\epsilon_1, \epsilon_2\}$ be a local orthonormal frame on 
$U \subset M$,  
and let $\nu$ denote a unit normal field. Regarding them as 
$\L^4$-valued functions, we consider a map 
\begin{equation*}
 (e_0,e_1,e_2,e_3):=(f,\epsilon_1,\epsilon_2, \nu) \colon 
U \to \mathcal F. 
\end{equation*} 
We shall use the same notation for differential forms on $\mathcal F$ 
and forms on $U$ pulled back by this map. 
Since $\la \nu, df \ra_L = 0$,   
\begin{equation}\label{eq:omega3=0}
 0 = \la \nu, df \ra_L 
=\la e_3, de_0 \ra_L 
= \omega^3.  
\end{equation}
From now on, we shall use the following convention 
on the ranges of indices: $1 \le i, j ,k \le 2$. 
It follows from \eqref{eq:conn_form-1}, \eqref{eq:de0ij}
and \eqref{eq:omega3=0} that 
\begin{align}
& \omega^{\alpha}_{\alpha}=0, \quad 
 -\omega_{i}^{0} + \omega_{0}^{i}=0, \quad
 -\omega_{3}^{0} + \omega_{0}^{3}=0, \quad 
 \omega_{i}^{j} + \omega^{i}_j =0, \quad
 \omega_{3}^{j} + \omega^{3}_j =0, \notag \\
& de_0 = e_{i} \otimes \omega^{i}, \quad
 de_i = e_{0} \otimes \omega^{i} + e_{j} \otimes \omega_{i}^{j}
        + e_{3} \otimes \omega^{3}_{i}, \quad
 de_3 =  e_{j} \otimes \omega_{3}^{j}.  \label{eq:de0i3}
\end{align}
And the structure equations \eqref{eq:str-1} become 
\begin{align}
& d \omega^i = -\omega^{i}_{j} \wedge \omega^j, \quad
0= d \omega^3 = -\omega^{3}_{j} \wedge \omega^j \label{eq:1st_str_eq}\\
& d \omega^1_2 = -\omega^1_3 \wedge \omega^3_2 - \omega^1 \wedge
 \omega^2, \quad
d \omega^3_j = -\omega^3_k \wedge \omega^k_j. \label{eq:2nd_str_eq} 
\end{align}
Following Bryant's notation (\cite{Br}), 
we introduce two complex-valued one-forms 
$\omega:=\omega^1 + \iii \omega^2$, 
$\pi := \omega^3_1 - \iii \omega^3_2$%
, and a complex vector $e := (e_1 - \iii e_2)/{2}$. Then   
\eqref{eq:1st_str_eq} and \eqref{eq:2nd_str_eq} are rewritten as 
\begin{align}
 d \omega &= \iii \omega^1_2 \wedge \omega, \quad 
\omega \wedge \pi + \bar \omega \wedge \bar \pi =0,  
 \label{al:domega} \\
 d \omega^1_2 &= - \frac{\iii}{2} 
(\pi \wedge \bar \pi + \omega \wedge \bar \omega), \quad 
d \pi = - \iii \omega^1_2 \wedge \pi.
 \label{al:dpi}
\end{align}

The first fundamental form $\fff = \la de_0, de_0 \ra_L$ is 
given by 
\begin{equation*}
\fff = \la e_{i} \otimes \omega^{i}, e_{i} \otimes \omega^{i} \ra_L
=(\omega^1)^2 + (\omega^2)^2 
= \omega \bar \omega = |\omega|^2. 
\end{equation*}
The Gaussian curvature $K$ is determined by 
$d \omega^1_2 = K \, ({\iii}/{2}) \, \omega \wedge \bar \omega $.  
Hence, it follows from \eqref{al:dpi} that 
\begin{equation}\label{eq:(K+1)}
 (K+1)\omega \wedge \bar \omega + \pi \wedge \bar \pi =0.
\end{equation}

The second fundamental form $\sff =  -\la de_0, de_3 \ra_L $ is 
given by 
\begin{equation*}
\sff =-\la e_{i} \otimes \omega^{i}, e_{i} \otimes \omega^{i}_3 \ra_L
=\omega^1\omega^3_1 + \omega^2\omega^3_2 
 = \frac{1}{2} (\omega \pi + \bar \omega \bar \pi)
= \text{Re}(\omega \pi).  
\end{equation*}
If we set $\omega^3_i = h_{ij}\omega^j$, then $h_{12}=h_{21}$ 
and $\sff =
h_{11}(\omega^1)^2 + 2 h_{12}\omega^1\omega^2 + h_{22}(\omega^2)^2$. 
Moreover, 
\begin{equation*}
  \pi = \frac{1}{2} \{ (h_{11}-h_{22}) -2 \iii h_{12} \} \omega 
      + \frac{1}{2}(h_{11}+h_{22}) \bar \omega.  
\end{equation*}
 Setting $q= \{ (h_{11}-h_{22}) -2i h_{12} \} /{2}$ and 
$H=(h_{11}+h_{22}) /{2}$, we can write 
\begin{equation}\label{eq:pi}
 \pi = q \omega + H \bar \omega. 
\end{equation}
Here, $H$ is the mean curvature. The second fundamental form $\sff$ 
is written as 
\begin{equation}\label{eq:sff}
 \sff = \frac{q}{2}\omega  \omega + H \omega \bar \omega
+ \frac{\bar q}{2}\bar \omega \bar \omega. 
\end{equation} 
It follows from \eqref{eq:(K+1)} and \eqref{eq:pi} that 
\begin{equation}\label{eq:K}
K=-1+H^2 - |q|^2 \left( = -1+\det (h_{ij})\right). 
\end{equation}
As a corollary, 
\begin{equation}\label{eq:H2-K-1}
 H^2 -K -1 \ge 0
\end{equation}
holds at every point $p \in M$, 
with equality if and only if $p$ is an umbilic point. 

The third fundamental form $\tff= \la de_3, de_3 \ra_L$ is given by 
\begin{equation*}
 \tff=\la e_{i} \otimes \omega^{i}_3, e_{i} \otimes \omega^{i}_3 \ra_L
=(\omega^1_3)^2 + (\omega^2_3)^2 
= \pi \bar \pi = |\pi|^2. 
\end{equation*}
\smallskip

The ideal boundary $\pd \H^3$ is considered as the quotient space 
$N^3 /{\sim} $, where  
\begin{equation*}
 N^3 = \setdef{x=(x_0,x_1,x_2,x_3)\in \L^4}{\la x, x \ra_L=0, \ x_0>0} 
\end{equation*}
and $x \sim y$ if $x = \lambda y$ for some positive constant 
$\lambda$. In other words, $\pd \H^3$ consists of
 positive null half-lines in $\L^4$.  
$N^3 /{\sim}$ is diffeomorphic to the $2$-sphere, and 
a natural conformal structure on $N^3 /{\sim}$ is given 
by the induced metric on $N^3$. 
Hence, $\pd \H^3 (= N^3 /{\sim})$ is identified 
with the conformal $2$-sphere.   
By definition, the {\it hyperbolic Gauss maps\/} are 
$G^{\pm}= [e_0 \pm e_3] \colon M \to \pd \H^3$, where 
$[v]$ denotes the line spanned by $v \in \L^4$. 
Because we mainly treat $G^+=[e_0+e_3]$, we simply write 
$G$ for $G^+$.  

The conformal structure on $M$ induced by $G$ is the conformal class 
determined by $\la d(e_0 + e_3), d(e_0 + e_3) \ra_L$. 
Indeed, it is computed as follows:
\begin{lemma}\label{eq_ffs}
\begin{equation*}\label{eq:ffs}
   \la d(e_0 + e_3), d(e_0 + e_3) \ra_L =|\omega -\bar \pi|^2 
= 2 (H - 1) \sff -K \fff.  
\end{equation*}
\end{lemma}
\begin{proof}
Since  
 $ de_0 + de_3 =e \otimes (\omega - \bar \pi)
+ \bar e \otimes (\bar \omega - \pi)$ holds by \eqref{eq:de0i3}, 
the first equality is obvious. The second equality follows from 
a straightforward computation using \eqref{eq:pi}, \eqref{eq:sff} 
and \eqref{eq:K}. 
\end{proof}
Similarly, the third fundamental form $\tff= \la de_3, de_3 \ra_L$
is computed as 
\begin{equation}\label{eq:tff}
 \tff= \la de_3, de_3 \ra_L = 2 H \sff -(K+1) \fff.
\end{equation}
\begin{proposition}\label{prop:tildeK}
 The Gaussian curvature $\tilde K$ of the 
pseudometric $\la d(e_0 +e_3), d(e_0 +e_3) \ra_L$  
is
\begin{equation*}
 \tilde K = \frac{K}{K-2(H-1)}. 
\end{equation*} 
\end{proposition}
\begin{proof}
 Setting $\alpha =\omega -\bar \pi$, we have  
\begin{equation*}
 d \alpha = d \omega - d \bar \pi 
= \iii \omega^1_2 \wedge \omega -(\iii \omega^1_2 \wedge \bar \pi)
= \iii \omega^1_2 \wedge \alpha. 
\end{equation*}
Hence, we can consider $\omega^1_2$ as the connection form of the 
metric $|\alpha|^2$. On the other hand, by \eqref{eq:pi} and 
\eqref{eq:K},  
\begin{align*}
  \alpha \wedge \bar \alpha 
&= (\omega -\bar \pi) \wedge (\bar \omega - \pi) 
= \{(1-H)\omega -\bar q \bar \omega)\} \wedge 
\{(1-H)\bar \omega - q \omega)\} \\
&= \{(1-H)^2 -|q|^2 \} \omega \wedge \bar \omega 
= (K-2H+2)\omega \wedge \bar \omega.  
\end{align*}
Therefore 
\begin{equation*}
  d \omega^1_2 = K \cdot \frac{\iii}{2}  \omega \wedge \bar \omega = 
 \frac{K}{K-2H+2}  \cdot \frac{\iii}{2} \alpha \wedge \bar \alpha, 
\end{equation*}
which proves the assertion. 
\end{proof}
\subsection{Parallel surfaces}
A map $f_t:=\cosh t \, f + \sinh t \, \nu$ is called the 
{\it parallel surface} of $f$ at distance $t$. 
It is easily verified that $f_t \colon M \to \H^3$ and that 
 $f_t(p)$ is joined to $f(p)$ by a hyperbolic line segment of 
length $t$. In general, $f_t$ may fail to be an immersion. In fact, 
$f_t$ is an immersion if and only if 
$\cosh t \, \omega^i + \sinh t \, \omega^i_3 \ne 0$ 
for every $p \in M$, 
because 
\begin{align*}
 df_t &= \cosh t \, df + \sinh t \, d \nu 
       = \cosh t \, de_0 + \sinh t \, de_3 \\
      &= \cosh t \, e_i \otimes \omega^i 
       + \sinh t \, e_i \otimes \omega^i_3 
      = e_i \otimes (\cosh t \, \omega^i + \sinh t \, \omega^i_3). 
\end{align*}

In this section, we assume that $f_t$ is an immersion,  
unless otherwise stated. 

\smallskip

The first fundamental form $\fff_t = \la df_t, df_t \ra_L$ is 
\begin{equation*}
 \fff_t = (\cosh t \ \omega^1 + \sinh t \ \omega^1_3)^2 
+ (\cosh t \ \omega^2 + \sinh t \ \omega^2_3)^2, 
\end{equation*}
hence, the $\theta^i := \cosh t \, \omega^i + \sinh t \, \omega^i_3$  
($i=1,2$) form an orthonormal frame of $f_t$. 
It follows from the structure equations \eqref{eq:1st_str_eq} 
and  \eqref{eq:2nd_str_eq} that 
\begin{align*}
 d \theta^i = -\omega^i_j \wedge \theta^j.  
\end{align*}
Thus $\omega^1_2$ is also a connection form of 
 $\fff_t = (\theta^1)^2 + (\theta^2)^2$. 
Denoting the Gaussian curvature of $\fff_t$ by $K_t$, we have 
\begin{equation*}
 d \omega^1_2 = K_t \theta^1 \wedge \theta^2 \ 
             (= K \omega^1 \wedge \omega^2). 
\end{equation*}
Using $K=-1 + \det(h_{ij})$ and $2H=h_{11}+h_{22}$, we have 
\begin{align*}
 \theta^1 \wedge \theta^2 
= \{ \cosh^2 t -2H \cosh t \sinh t  
+ (K+1)\sinh^2 t  \} \,
\omega^1 \wedge \omega^2.  
\end{align*}
Therefore 
\begin{align*}
 d \omega^1_2 = \frac{K}{K\sinh^2 t-2H \cosh t \sinh t  
+\cosh^2 t + \sinh^2 t } \theta^1 \wedge \theta^2. 
\end{align*} 
This implies that 
\begin{align}\label{K_t}
 K_t = \frac{K}{K\sinh^2 t-2H \cosh t \sinh t  
+\cosh^2 t + \sinh^2 t }. 
\end{align}
Thus, since $(K_t)_{-t}=K$, we have  
 \begin{align*}
 \frac{K_t}{K_t \sinh^2 t + 2H_t \cosh t \sinh t  
+\cosh^2 t + \sinh^2 t } 
= K. 
\end{align*}
This formula together with \eqref{K_t} implies that 
\begin{align}\label{H_t}
 H_t = \frac{H(\cosh^2 t + \sinh^2 t)-(2+K) \cosh t \sinh t}
{K\sinh^2 t-2H \cosh t \sinh t +\cosh^2 t + \sinh^2 t }. 
\end{align}
The formulas \eqref{K_t} and \eqref{H_t} yield the following 
well-known theorem:
\begin{theorem}\label{th:para-1}
 \begin{enumerate}
  \item[\rm (1)] All parallel surfaces of a flat surface are also flat. 
  \item[\rm (2)] A family of parallel surfaces of a 
surface with constant mean curvature $(|H|>1)$ contains 
a surface with constant Gaussian curvature $(K>0)$,  and vice versa. 
 \end{enumerate}
\end{theorem}

\smallskip

We can rewrite \eqref{H_t} as 
\begin{equation*}
 H_t-1 = \frac{(\cosh t+\sinh t)\{(\cosh t+\sinh t)(H-1) -K \sinh t\}}
{K\sinh^2 t-2(H-1) \cosh t \sinh t +(\cosh t - \sinh t)^2 }. 
\end{equation*}
Multiplying $K$ on both sides, we have 
\begin{equation*}
 K(H_t-1) = K_t \, e^t\{e^t(H-1) -K \sinh t\}. 
\end{equation*}
For example, if we assume that the original surface $f$ has constant 
mean curvature one (CMC-1), then 
\begin{equation*}
 H_t-1 = (-e^t \sinh t) K_t,  
\end{equation*}
thus, $f_t$ is a Weingarten surface. The family of 
Weingarten surfaces satisfying $H-1 = \lambda K$ for some constant 
$\lambda$ includes the following interesting surfaces:
\begin{quote}
 If $\lambda=0$, then $f$ is a  CMC-1 surface. \\
 If $\lambda=1/2$, then at least one of the principal curvatures 
equals $1$. \\
 If $\lambda=1$, then the sum of the reciprocals of the principal 
curvature is the constant value $2$, that is, $f$ is a
surface with constant harmonic-mean curvature one (HMC-1).   
\end{quote}
These are verified by $K = -1 + \kappa_1 \kappa_2$,  
$2H = \kappa_1 + \kappa_2$ and that the harmonic-mean curvature is 
$2/(\kappa_1^{-1} + \kappa_2^{-1})$,
 where $\kappa_i$ ($i=1,2$) denote 
the principal curvatures.

Conversely, we assume that the original surface $f$ satisfies 
$H-1 = \lambda K$ for some constant $\lambda$. 
Then the parallel surface $f_t$ satisfies
\begin{equation*}
 H_t -1 = \frac{(2 \lambda -1)e^{2t} +1}{2}K_t, 
\end{equation*}
and hence is the same kind of Weingarten surface. 
Since $\lambda_t := \{(2 \lambda -1)e^{2t} +1\}/{2}$ satisfies 
$(2 \lambda_t -1)=(2 \lambda -1)e^{2t}$, the following lemma is  
clear:
\begin{lemma}
 \begin{enumerate}[\rm (1)]
 \item If $\lambda=1/2$, then $\lambda_t = 1/2$ for all $t$. 
 \item If $\lambda<1/2$, then $\lambda_t < 1/2$ for all $t$, and
 $\lambda_t =0$ for some unique $t$. 
 \item If $\lambda>1/2$, then $\lambda_t > 1/2$ for all $t$, and  
 $\lambda_t =1$ for some unique $t$. 
\end{enumerate}
\end{lemma}
Therefore, we have the following theorem: 
\begin{theorem}\label{th:para-2}
 \begin{enumerate}[\rm (1)]
 \item  Let $f$ be a surface 
satisfying that at least one of the principal curvatures 
equals $1$, i.e., a Weingarten surface with $H-1 = {K}/{2}$. 
Then all parallel surfaces of $f$ also  
satisfy that at least one of the principal curvatures 
equals $1$. 
 \item Let $f$ be a Weingarten surface with $H-1 = \lambda K$ for 
some constant $\lambda (>1/2)$. 
Then the family of parallel surfaces of $f$ consists of 
Weingarten surfaces with $H-1 = \lambda K \ (\lambda >1/2)$.  
This family includes a single HMC-1 surface.  
 \item Let $f$ be a Weingarten surface with $H-1 = \lambda K$ for 
some constant $\lambda (<1/2)$. 
Then the family of parallel surfaces of $f$ consists of 
Weingarten surfaces with $H-1 = \lambda K \ (\lambda <1/2)$.  
This family includes a single CMC-1 surface.  
 \end{enumerate}
\end{theorem}
 Theorem \ref{th:para-1} is well-known, 
whereas Theorem \ref{th:para-2} seems to be lesser-known. 
\subsection{Weingarten surfaces satisfying $H-1=\lambda K$}
Throughout this section, $f \colon M \to \H^3$ denotes a Weingarten surface 
satisfying $H-1=\lambda K$ for some constant $\lambda$, unless otherwise
stated.

\medskip

It follows from Lemma \ref{eq_ffs} that  
\begin{equation*}
   \la d(e_0 + e_3), d(e_0 + e_3) \ra_L = 
2 (H - 1) \sff -K \fff 
= -K(\fff - 2 \lambda \sff ). 
\end{equation*}
Hence, if we endow $M$ with the ``metric'' $\fff - 2 \lambda \sff$,  
then the hyperbolic Gauss map $G$ is conformal. However, we need to
check that $\fff - 2 \lambda \sff$ is indeed a metric:  
\begin{lemma}
If $\lambda \ne 1/2$, then 
$\fff - 2 \lambda \sff$ is either positive or negative definite. 
\end{lemma}
\begin{proof}
$\fff - 2 \lambda \sff $ is definite if and only if 
\begin{equation}\label{eq:det1-22}
 \det \begin{bmatrix}
      1-2 \lambda h_{11} & -2 \lambda h_{12} \\ 
      -2 \lambda h_{12} & 1-2 \lambda h_{22}
     \end{bmatrix} >0.  
\end{equation}
This condition \eqref{eq:det1-22} is equivalent to
$4 \lambda^2 ( K+1) -2 \lambda ( 2H) +1 > 0$. Moreover, from 
the assumption $H-1 = \lambda K$, this is equivalent to 
$(2 \lambda -1)^2 >0$. 
\end{proof}
 
As stated before, the special case $\lambda =0$ concerns  
CMC-1 surfaces. 
The following proposition is shown in \cite{Br} when $\lambda=0$, and 
can be proved by the same argument. Thus, the proof is omitted here. 
\begin{proposition}
Let $f \colon M \to \H^3$ be a Weingarten surface 
satisfying $H-1=\lambda K$ for some constant $\lambda (\ne 1/2)$. 
Then $\fff -2 \lambda \sff$ determines a conformal structure on $M$, and 
the hyperbolic Gauss map $G \colon (M, \fff -2 \lambda \sff) \to \pd \H^3$ 
is conformal.

Conversely, if an immersed surface $f \colon M \to \H^3$ satisfies 
\begin{enumerate}[\rm (i)]
 \item $\fff - 2 \lambda_0 \sff$ is definite, and 
 \item $G \colon (M, \fff -2 \lambda_0 \sff ) \to \pd \H^3$ 
is conformal 
\end{enumerate}   
for some constant $\lambda_0$, then 
$f$ is a totally umbilic surface or a Weingarten surface 
satisfying $H-1=\lambda_0 K$.  
\end{proposition} 

The following propositions follow easily from \eqref{eq:H2-K-1}
and Proposition \ref{prop:tildeK}, respectively.
\begin{proposition}\label{prop:KH_BLW}
If $f \colon M \to \H^3$ is a Weingarten surface satisfying 
$H-1=\lambda K$ $(\lambda \ne 1/2)$, then the Gaussian curvature  
$K$ satisfies the following inequalities: 
\begin{enumerate}[\rm (i)]
 \item If $\lambda =0$, then $K\le 0$. 
 \item If $\lambda <1/2 (\ne 0)$, then $K\le 0$ or 
$K \ge (1 -2 \lambda)/\lambda^2$. 
 \item If $\lambda >1/2$, then  
$K \le (1 -2 \lambda)/\lambda^2$ or $K\ge 0$. 
\end{enumerate}
And the mean curvature $H$ satisfies the following: 
\begin{enumerate}[\rm (i)]
 \item If $0 < \lambda <1/2$, then $H\le 1$ or 
$H \ge (1 - \lambda)/\lambda$. 
 \item If $\lambda < 0$ or $\lambda >1/2$, then  
$H \le (1 - \lambda)/\lambda$ or $H \ge 1$. 
\end{enumerate}
\end{proposition}
\begin{proposition}
For a Weingarten surface satisfying 
$H-1=\lambda K$ $(\lambda \ne 1/2)$,  
the pseudometric $\la d(e_0+e_3),d(e_0+e_3) \ra_L $ 
has constant curvature $1/(1-2 \lambda)$.
\end{proposition} 
%
%
\section{Surfaces with constant harmonic-mean curvature one}\label{sec:mrc-1surf}
\subsection{Basics}
In this section, we study the case $\lambda=1$ for 
Weingarten surfaces satisfying $H-1=\lambda K$, that is, 
the case $H-1=K$.   
As stated in the previous section, a Weingarten surface satisfying 
$H-1=K$ has constant harmonic-mean curvature one, 
and we call it an HMC-1 surface.

By Lemma \ref{eq_ffs} and \eqref{eq:tff}, an HMC-1 surface satisfies
\begin{align}
 \la de_0 + de_3, de_0 + de_3 \ra_L 
&= -K(\fff -2 \sff), \label{al:de_0+de_3de_0+de_3}\\ 
(\tff=) \la de_3, de_3 \ra_L &= -H(\fff -2 \sff). \label{al:de_3de_3}
\end{align}
These two quadratic differentials are conformally equivalent. 
Following Bryant's notation (\cite{Br}), we set 
\begin{equation*}
 \eta := (\omega^1 - \omega^3_1) - \iii (\omega^2 - \omega^3_2)
\ (= \bar \omega - \pi). 
\end{equation*}
The formulas \eqref{al:de_0+de_3de_0+de_3} and \eqref{al:de_3de_3}
become 
\begin{align}
|\eta |^2 & \left( = 
 (\omega^1 -\omega^3_1)^2 + (\omega^2 -\omega^3_2)^2 \right) = 
-K(\fff -2 \sff), \label{al:|eta|}\\ 
|\pi|^2 & \left( = (\omega^3_1)^2 + (\omega^3_2)^2 \right) = 
-H(\fff -2 \sff). \label{al:|pi|}
\end{align}
Consequently, we obtain the following lemma:
\begin{lemma}
  $\eta(p)=0$ if and only if $p$ is an umbilical point 
with $\sff_p = \fff_p$, i.e., a point 
where $K=0$ and $H=1$.  
  $\pi(p)=0$ if and only if $p$ is a totally geodesic point, 
i.e., a point where $K=-1$ and $H=0$. 
\end{lemma}
\begin{lemma}\label{lem:KH} 
The Gaussian and mean curvatures $K$, $H$ are given by 
\begin{equation*}
 K = \frac{|\eta|^2}{|\pi|^2 - |\eta|^2}, \quad 
 H = \frac{|\pi|^2}{|\pi|^2 - |\eta|^2}. 
\end{equation*}
\end{lemma}
\begin{proof}
It follows from \eqref{al:|eta|} and \eqref{al:|pi|} that 
 $H |\eta|^2 = K |\pi|^2$. This formula implies the assertion, 
since $H-1=K$. 
\end{proof}
Note that 
\begin{equation*}
 K \le -1 \ (H \le 0) \text{ or } 
 K \ge 0 \ (H \ge 1) 
\end{equation*}
for HMC-1 surfaces, because of Proposition \ref{prop:KH_BLW}. 
It follows from Lemma \ref{lem:KH} that  
$K \le -1$ if and only if $|\pi|^2 < |\eta|^2$, and that   
$K \ge 0$ if and only if $|\pi|^2 > |\eta|^2$.
\begin{lemma}
 \begin{equation}\label{eq:barq-eta}
  \bar q \eta = -K \bar \pi .
 \end{equation}
\end{lemma}
\begin{proof}
By \eqref{eq:K}, we have
$|q|^2 = H^2 -(K +1) = H^2 -H = H(H-1) = HK$.
Moreover, using \eqref{eq:pi}, we have 
$ \eta = \bar \omega - \pi = \bar \omega - (q \omega + H \bar \omega) 
=-q \omega + (1-H) \bar \omega = -q \omega -K \bar \omega$, 
therefore, 
$ \bar q \eta = -|q|^2 \omega -K \bar q \bar \omega 
= -HK \omega -K \bar q \bar \omega =(-K)(H \omega + \bar q \bar \omega )
=(-K) \bar \pi$.  
\end{proof}

\medskip

Now we shall equip $M$ with the complex structure that is compatible 
with the conformal structure $[\fff -2 \sff]$. 
Here, we give the orientation of $M$ so that  
$\eta$ is a $(1,0)$-form. Note that $\bar \pi$ is also a $(1,0)$-form, 
because of \eqref{eq:barq-eta}. 
\begin{lemma}\label{lem:Hopf} 
Let $\fff^{2,0}$ and $\sff^{2,0}$ 
denote the $(2,0)$-parts of the complexification of 
the fundamental forms $\fff$ and $\sff$, respectively. Then  
\begin{equation}\label{eq:fff20}
 \fff^{2,0} = {2} \sff^{2,0} = \eta \bar \pi, 
\end{equation}
and $\eta \bar \pi$ is a holomorphic quadratic differential on $M$. 
\end{lemma}
\begin{proof}
\eqref{eq:fff20} is obtained by rewriting $\fff$ and  
$\sff$ with $\eta = \bar \omega - \pi$. Indeed, 
\begin{equation*}
 \fff \, (= \omega \bar \omega) = \eta \bar \pi + |\eta|^2 + |\pi|^2 
+ \bar \eta \pi, \quad 
\sff \, (=\text{Re}(\omega \pi)) = 
\frac{1}{2}(\eta \bar \pi + 2 |\pi|^2 
+ \bar \eta \pi). 
\end{equation*}
It follows from the formulas in \eqref{al:domega}, \eqref{al:dpi} that 
$d \bar \pi = \iii \omega^1_2 \wedge \bar \pi$ and 
$d \eta = - \iii \omega^1_2 \wedge \eta$. 
They imply that $\eta \bar \pi$ is holomorphic. 
\end{proof}
\begin{proposition}
 The pseudometric $|\pi|^2(= \!\tff)$ has the Gaussian curvature \\
 ${K}/(K+1) \, (=|\eta|^2 / |\pi|^2)$. 
\end{proposition}
\begin{proof}
Since $d \bar \pi = \iii \omega^1_2 \wedge \bar \pi$, 
we may regard $\omega^1_2$ as the connection form for 
$|\pi|^2$. Moreover, its exterior differential is computed as
\begin{align*}
 d \omega^1_2 &= - \frac{\iii}{2} 
(\pi \wedge \bar \pi + \omega \wedge \bar \omega)  
= - \frac{\iii}{2} 
(\pi \wedge \bar \pi + \frac{1}{K+1} \bar \pi \wedge \pi) \\
&= - \frac{\iii}{2} 
(\frac{-K}{K+1} \bar \pi \wedge \pi)
=\frac{K}{K+1} \frac{\iii}{2} \bar \pi \wedge \pi, 
\end{align*} 
because of \eqref{al:dpi} and \eqref{eq:(K+1)}, 
which proves the assertion.
\end{proof}
\begin{lemma}\label{lem:barpdpd}
 \begin{equation*}
  \bar \pd \pd (e_0+e_3) =
   \frac{1}{2}e_0 \otimes (\bar \eta \wedge \eta ).  
 \end{equation*}
\end{lemma}
\begin{proof}
Because $d(e_0+e_3) = e \, \bar \eta + \bar e \, \eta$, we have 
$\pd (e_0+e_3) = \bar e \, \eta$. Furthermore, taking $\bar \pd$, 
we can calculate as follows: 
 \begin{align*}
\bar \pd \pd (e_0+e_3) & = 
\bar \pd \bar e \wedge \eta + \bar e \otimes \bar \pd \eta \\ 
&= d \bar e \wedge \eta + \bar e \otimes d \eta \quad 
(\text{since $\eta$ is a $(1,0)$-form,})\\ 
&= \left(e_0 \otimes \frac{1}{2} \omega +\iii \bar e \otimes \omega^1_2 
+e_3 \otimes \frac{1}{2} \bar \pi \right) 
\wedge \eta + \bar e \otimes \left(-\iii \omega^1_2 \wedge \eta \right) \\
&= \frac{1}{2}e_0 \otimes (\omega \wedge \eta ) \quad
(\text{since $\bar \pi$ is a $(1,0)$-form,}) \\ 
&= \frac{1}{2}e_0 \otimes ((\omega - \bar \pi) \wedge \eta ) \quad 
(\text{since $\bar \pi$ is a $(1,0)$-form,}) \\
&= \frac{1}{2}e_0 \otimes (\bar \eta \wedge \eta ). 
\end{align*}
\end{proof}
\subsection{An overview on the work 
by G\'alvez, Mart\'\i{}nez and Mil\'an}\label{sec:ov-GMM}
In this section we give an overview on 
deriving G\'alvez-Mart\'\i{}nez-Mil\'an's formula 
stated in the introduction, restricting ourselves to 
HMC-1 surfaces. 

\medskip

First of all, we review the matrix model for $\H^3$. 
We identify $\L^4$ with $\Herm(2)$, the set of $2\times 2$ Hermitian
matrices, via 
\[
      \L^4\ni x=(x_0,x_1,x_2,x_3)\leftrightarrow
      X=\begin{bmatrix}
        x_0+x_3 & x_1+\sqrt{-1}x_2 \\
        x_1-\sqrt{-1}x_2 & x_0-x_3
      \end{bmatrix}\in\Herm(2).
\]
Since $\la x, x \ra_L = -\det X$ and $2 x_0 = \trace X$, 
\begin{align*}
     \H^3 &= \{ X\in\Herm(2)\,;\,\det X=1,~\trace X>0\}\\
         &= \{ aa^*\,;\,a\in\SL(2,\C)\}=\SL(2,\C)/\SU(2), \\
    \pd \H^3 &= \{ X\in\Herm(2)\,;\,\det X=0,~\trace X>0\} /{\sim} \\
         &= \left\{ aa^*\,;\,a=\begin{bmatrix} a_1 \\ a_2 \end{bmatrix} 
\in \C^2 \setminus \{ 0 \} \right\} /{\sim} \\ 
&= (\C^2 \setminus \{ 0 \} )/(\C \setminus \{ 0 \} )=\C P^1, 
\end{align*}
where $a^*$ is the conjugate transpose of $a$, 
and $\C P^1$ is the complex projective line. 

Hereafter, we will consider $\H^3$ to be $\SL(2,\C)/\SU(2)$, 
and $\pd \H^3$ to be $\C P^1$.   

\smallskip

Recall that $G=[e_0+e_3] \colon M \to \C P^1(= \pd \H^3)$ is 
a conformal map. 
Thus 
\begin{equation}\label{eq:e0+e3}
 e_0+e_3 = \Lambda \begin{bmatrix}
		     A \\ B 
		   \end{bmatrix}
\begin{bmatrix}
\, \bar A & \bar B \, 
\end{bmatrix}
\end{equation}
for some holomorphic functions $A$, $B$ and a positive 
function $\Lambda$. Note that 
$A$, $B$ and $\Lambda$ have an ambiguity, but  
$\Lambda |A|^2$, $\Lambda |B|^2$ and $\Lambda A \bar B$ 
are well-defined on $M$.  

Let $z$ be a local holomorphic coordinate on $M$, and 
let the lower suffix 
denote the partial derivative. 
Differentiating \eqref{eq:e0+e3} twice, we have 
\begin{equation*}
 (e_0+e_3)_{z \bar z}
=\begin{bmatrix}
  A & A_z \\ B & B_z
 \end{bmatrix}
\begin{bmatrix}
 \Lambda_{z \bar z} & \Lambda_{z} \\
 \Lambda_{\bar z}   & \Lambda
\end{bmatrix}
\begin{bmatrix}
  \lbar A & \lbar B \\ 
\lbar{A_z} & \lbar{B_z}
 \end{bmatrix}.
\end{equation*}
Hence, this and Lemma \ref{lem:barpdpd} imply that  
\begin{equation*}
 e_0 = \frac{2}{|\eta/dz|^2}
\begin{bmatrix}
  A & A_z \\ B & B_z
 \end{bmatrix}
\begin{bmatrix}
 \Lambda_{z \bar z} & \Lambda_{z} \\
 \Lambda_{\bar z}   & \Lambda
\end{bmatrix}
\begin{bmatrix}
  \lbar A & \lbar B \\ 
\lbar{A_z} & \lbar{B_z}
 \end{bmatrix}. 
\end{equation*}
Here we have assumed (and will continue to assume) that 
$\eta$ is not identically zero. This assumption means that the surface is 
not a horosphere. 

Setting 
\begin{equation*}
 g :=\begin{bmatrix}
  A & A_z \\ B & B_z
 \end{bmatrix}, \quad \delta := \frac{2}{|\eta/dz|^2}, \quad 
 \Omega := \delta \begin{bmatrix}
 \Lambda_{z \bar z} & \Lambda_{z} \\
 \Lambda_{\bar z}   & \Lambda
\end{bmatrix}, 
\end{equation*}
we have the following formulas:  
\begin{equation}\label{eq:e_0:1}
 e_0 = g \Omega g^{*}, 
\end{equation}
\begin{equation}\label{eq:e_3:1}
 e_3 \ (=(e_0+e_3)-e_0)=g \begin{bmatrix}
			   \Lambda & 0 \\ 0 & 0
			  \end{bmatrix} g^{*}
-g \Omega g^{*}
=g \tilde \Omega g^{*}, 
\end{equation}
where
\begin{equation*}
 \tilde \Omega = \begin{bmatrix}
			   \Lambda & 0 \\ 0 & 0
			  \end{bmatrix} - \Omega
= \begin{bmatrix}
\Lambda - \delta \Lambda_{z \bar z} & - \delta \Lambda_{z} \\
-\delta \Lambda_{\bar z}   & -\delta \Lambda
\end{bmatrix}.
\end{equation*}
\begin{lemma}\label{lem:eta-lambda}
 \begin{equation*}
  |\eta/dz|^2 = | \det g |^2 \Lambda^2 = 4 (\log \Lambda)_{z \bar z} .
 \end{equation*}
\end{lemma}
\begin{proof}
 \begin{align*}
  -1 &= \la e_0, e_0 \ra_L = - \det e_0 = - \det (g \Omega g^{*}) 
= -|\det g|^2 \delta^2 (\Lambda_{z \bar z}\Lambda - \Lambda_z
  \Lambda_{\bar z}), \\
   1 &= \la e_3, e_3 \ra_L = - \det e_3 = - \det (g \tilde \Omega g^{*})
= -|\det g|^2 \left\{ -\delta \Lambda^2 
+ \delta^2 (\Lambda_{z \bar z}\Lambda - \Lambda_z \Lambda_{\bar z}) \right\}.
 \end{align*}
Subtracting and adding these, we have 
\begin{align*}
 -2 &= -|\det g|^2 \delta \Lambda^2,  \text{ i.e., } 
|\eta/dz|^2 = | \det g |^2 \Lambda^2, \\
 0 &= -|\det g|^2 \left\{ -\delta \Lambda^2 
+2 \delta^2 (\Lambda_{z \bar z}\Lambda - \Lambda_z \Lambda_{\bar z})
\right\}, 
 \text{ i.e., } \Lambda^2 =
2\delta (\Lambda_{z \bar z}\Lambda - \Lambda_z \Lambda_{\bar z}).
\end{align*}
\end{proof}

By Proposition \ref{prop:tildeK}, the pseudometric 
$|\eta|^2=\la d(e_0+e_3), d(e_0+e_3) \ra_L$ has constant curvature $-1$ 
for an HMC-1 surface.  
It follows from the Frobenius theorem that there exists a holomorphic 
map $h$ from the universal cover $\tilde M$ to the Poincar\'e disk 
$\D$ such that the pull-back of the Poincar\'e metric via $h$ coincides 
with $|\eta|^2$, that is,  
\begin{align}
  1-|h|^2 > 0, \quad 
  |\eta|^2 = \frac{4 |dh|^2}{(1-|h|^2)^2}. \label{al:etaetabar}
\end{align}
  
Exchanging with $h$, we reexamine \eqref{eq:e_0:1} and \eqref{eq:e_3:1}. 
Since the pair $A$, $B$ has the ambiguity of multiplication by 
non-zero holomorphic functions, we can start with the assumption 
 \begin{equation}\label{eq:ABz-AzB}
 (\det g =) AB_z - A_z B = h_z. 
\end{equation}
It follows from Lemma \ref{lem:eta-lambda} and \eqref{al:etaetabar} that
\begin{equation}\label{eq:lambda}
 \Lambda = \frac{2}{1-|h|^2}. 
\end{equation}

It is straightforward to calculate that  
\begin{align}
 \Lambda_z &= \frac{2 h_z \bar h}{(1-|h|^2)^2}, \label{al:lambda_z} \\ 
 \Lambda_{z \bar z} &= \frac{2 |h_z|^2 (1+|h|^2)}{(1-|h|^2)^3}.  
\label{al:lambda_zzbar}
\end{align}
Substituting 
\eqref{al:etaetabar}, \eqref{eq:lambda}, \eqref{al:lambda_z}, 
\eqref{al:lambda_zzbar} into \eqref{eq:e_0:1}, \eqref{eq:e_3:1}, 
respectively, we have
\begin{align*}
 e_0 & = \begin{bmatrix}
  A & A_z \\ B & B_z
 \end{bmatrix}
\begin{bmatrix}
 \frac{1+|h|^2}{1-|h|^2} & \lbar{{h}/{h_z}} \\
 {h}/{h_z} & \frac{1-|h|^2}{|h_z|^2}
\end{bmatrix}
\begin{bmatrix}
  \lbar A & \lbar B \\ 
\lbar{A_z} & \lbar{B_z}
 \end{bmatrix} \\ 
&=\begin{bmatrix}
  A & A_z/h_z \\ B & B_z/h_z
 \end{bmatrix}
\begin{bmatrix}
 \frac{1+|h|^2}{1-|h|^2} & {\bar h} \\
 {h} & {1-|h|^2}
\end{bmatrix}
\begin{bmatrix}
  \lbar A & \lbar B \\ 
\lbar{A_z}/\lbar{h_z} & \lbar{B_z}/\lbar{h_z}
 \end{bmatrix}, \\
e_3 &= 
\begin{bmatrix}
  A & A_z/h_z \\ B & B_z/h_z
 \end{bmatrix}
\begin{bmatrix}
 1 & -{\bar h} \\
 -{h} & {-1+|h|^2}
\end{bmatrix}
\begin{bmatrix}
  \lbar A & \lbar B \\ 
\lbar{A_z}/\lbar{h_z} & \lbar{B_z}/\lbar{h_z}
 \end{bmatrix}. 
\end{align*}
Introducing the three matrices 
\begin{equation*}
 \G:=\begin{bmatrix}
  A & A_z/h_z \\ B & B_z/h_z
 \end{bmatrix}, \quad 
\HH:=\begin{bmatrix}
 \frac{1+|h|^2}{1-|h|^2} & {\bar h} \\
 {h} & {1-|h|^2}
\end{bmatrix}, \quad 
\tilde \HH :=
\begin{bmatrix}
 1 & -{\bar h} \\
 -{h} & {-1+|h|^2}
\end{bmatrix}, 
\end{equation*}
we can write 
\begin{equation*}
 e_0 = \G \HH \G^*, \quad e_3 = \G \tilde \HH \G^*. 
\end{equation*}
 By straightforward calculation using \eqref{eq:ABz-AzB}, i.e.,  
$A_z B -A B_z = h_z$, we have  
\begin{equation*}
  \G^{-1} d\G = \begin{bmatrix}
		0 & \theta \\ dh & 0
	       \end{bmatrix}, 
\text{ where }
\theta = \frac{B_z A_{zz}-A_z B_{zz}}{(h_z)^2}dz. 
\end{equation*}
The one-form $\theta$ is also written as 
\begin{equation}\label{eq:theta}
 \theta = \frac{1}{A} \, d \left( \frac{dA}{dh} \right)
= \frac{1}{B} \, d \left( \frac{dB}{dh} \right). 
\end{equation}
Note that $\theta$ is a one-form defined on $\tilde M$. 

\smallskip

In the following, we describe the fundamental forms in terms of 
$h$ and $\theta$. 
It is not difficult to calculate that 
\begin{align*}
 d e_0 =\G  
\begin{bmatrix}
 * & {2d \bar h}/{(1-|h|^2)} + (1-|h|^2)\theta  \\ 
{2dh}/{(1-|h|^2)}  +(1-|h|^2) \bar \theta & 0  
\end{bmatrix}
\G^*. 
\end{align*}
Hence, we have
\begin{align}\label{eq:fff}
 \fff = - \det(de_0) &= \left| 
\frac{2}{1-|h|^2} d h +(1-|h|^2) \bar \theta \right|^2 \\
&= \frac{4|dh|^2}{(1-|h|^2)^2} + 2 \theta dh 
+ 2 \bar \theta d \bar h +(1-|h|^2)^2 |\theta|^2. \notag 
\end{align}
As a by-product of this formula, we obtain the following lemma:
\begin{lemma}
 $\theta dh$ and $(1-|h|^2)^2|\theta|^2$ are well-defined on $M$.  
\end{lemma}
It is not difficult to calculate that 
\begin{align*}
 d e_3 
=\G 
\begin{bmatrix}
 * & (|h|^2-1)\theta \\ 
(|h|^2-1) \bar \theta & 0  
\end{bmatrix}
\G^*.  
\end{align*}
Hence, we have
\begin{equation*}
 \tff = - \det(de_3)= (1-|h|^2)^2 \left| \theta \right|^2. 
\end{equation*}
It follows that  
\begin{align}
 \fff - 2 \sff &= |de_0|_L^2 + \left\{
|de_0 + de_3 |_L^2 -|de_0|_L^2 - |de_3 |_L^2
\right\}=|de_0 + de_3 |_L^2 -|de_3|_L^2 \notag \\
&= \frac{4 |dh|^2}{(1-|h|^2)^2}-(1-|h|^2)^2|\theta|^2. \label{eq:fff-2sff}
\end{align}

From the argument above, 
one can understand the G\'alvez-Mart\'\i{}nez-Mil\'an formula 
stated in the introduction.

\bigskip

We finish this section by providing some other formulas.  \\
By \eqref{eq:fff} and \eqref{eq:fff-2sff}, we have
\begin{equation*}
 \sff = \theta dh +(1-|h|^2)^2 |\theta|^2 + \bar \theta d \bar h. 
\end{equation*}
In particular, 
\begin{equation}\label{eq:|pi|}
|\pi|^2=\sff^{1,1} = (1-|h|^2)^2 |\theta|^2. 
\end{equation}
It follows from Lemma \ref{lem:KH}, 
\eqref{al:etaetabar} and \eqref{eq:|pi|} that 
\begin{align*}
  K = \frac{4|dh|^2}{(1-|h|^2)^4 |\theta|^2-4|dh|^2}, \quad 
 H = \frac{(1-|h|^2)^4 |\theta|^2}{(1-|h|^2)^4 |\theta|^2-4|dh|^2}. 
\end{align*}
\subsection{Improvement of the representation formula}
We shall give a slight improvement of the  
G\'alvez-Mart\'\i{}nez-Mil\'an formula, limiting ourselves to
HMC-1 surfaces, 
and make it clear what is a local invariant and 
what is a global invariant. 

The hyperbolic Gauss map $G=A/B$ is globally-defined on $M$. 
We can represent $\G$ using $G$ as follows:
\begin{lemma}\label{lem:mathcalG}
 \begin{equation}\label{eq:mathcalG}
  \G = \left( -G_h \right)^{-{3}/{2}} \begin{bmatrix}
-G G_h & GG_{hh}/2 - G_h^2 \\
-G_h &  G_{hh}/2	
       \end{bmatrix}, 
 \end{equation}
 where $ G_h = {dG}/{dh}, \ G_{hh} = {d^2G}/{dh^2}$.

\end{lemma} 
\begin{proof} $\G$ is computed as 
 \begin{equation}\label{eq:mtrxG}
  \G = \begin{bmatrix}
      A & dA/dh \\ B & dB / dh
       \end{bmatrix}
     = \begin{bmatrix}
      GB & d(GB)/dh \\ B & dB / dh
       \end{bmatrix}
     = \frac{1}{B} \begin{bmatrix}
      GB^2 & Bd(GB)/dh \\ B^2 & BdB / dh
                    \end{bmatrix}.
 \end{equation}
On the other hand, substituting $A=BG$ to $AdB -B dA = dh$, we have
\begin{equation}\label{eq:B^2}
 B^2=- \dfrac{dh}{dG}=-\dfrac{1}{G_h}. 
\end{equation}
 Eliminating $B$ from 
\eqref{eq:mtrxG} with this, we have the assertion.  
\end{proof}
Thus, for an HMC-1 surface $f$, 
we can make a representation formula $f=\G \HH \G^*$ 
with $\G$ as in \eqref{eq:mathcalG} and 
\begin{equation}\label{eq:mtrixH}
 \HH=\begin{bmatrix}
 \frac{1+|h|^2}{1-|h|^2} & {\bar h} \\
 {h} & {1-|h|^2}
\end{bmatrix}
\end{equation}
from a meromorphic function $G$ on $M$ and a holomorphic map 
$h\colon \tilde M \to \D$. However, it is not defined on $M$ yet 
(merely on $\tilde M$, in general).  
We need to find the condition that $f$ is single-valued on $M$. 
Indeed, we prove:
\begin{proposition}
$f=\G \HH \G^*$ with $\G$, $\HH$ as in \eqref{eq:mathcalG}, \eqref{eq:mtrixH}
is single-valued on $M$ if and only if 
the pseudometric $|\eta|^2=4|dh|^2/(1-|h|^2)^2$ is 
single-valued on $M$. 
\end{proposition}
\begin{proof}
Suppose that $|\eta|^2=4|dh|^2/(1-|h|^2)^2$ is single-valued on $M$. 

By \eqref{eq:lambda} and Lemma \ref{lem:mathcalG}, we have 
\begin{align*}
 \Lambda |A|^2 &= \frac{2}{1-|h|^2}|G|^2|-G_h|^{-1} 
                = \frac{2}{1-|h|^2}|G|^2|\frac{dh}{dG}| 
                = 2 \frac{|dh|}{1-|h|^2} \frac{|G|^2}{|dG|} \\
 \Lambda |B|^2 &= 2 \frac{|dh|}{1-|h|^2} \frac{1}{|dG|} \\
 \Lambda A \bar B &= 2 \frac{|dh|}{1-|h|^2} \frac{G}{|dG|}. 
\end{align*}
Thus, all 
$\lambda |A|^2$, $\lambda |B|^2$, $\lambda A \bar B$ 
are single-valued on $M$. 
In other words, $e_0 + e_3$ is single-valued on $M$ 
(because of \eqref{eq:e0+e3}). 
On the other hand, $\eta \wedge \bar \eta$ 
is also single-valued on $M$. 
Therefore, recalling the formula  
$\bar \pd \pd (e_0+e_3)=\frac{1}{2}e_0 \otimes \bar \eta \wedge \eta$, 
we can conclude that $e_0 (=f)$ is also single-valued on $M$. 
\end{proof}

Therefore we have:
\begin{theorem}\label{thm:main1}
 Let $G$ be a meromorphic function on a Riemann surface $M$, 
and $|\eta|^2$ a pseudometric on $M$ of constant curvature $-1$. 
Suppose that the quadratic differential form \eqref{eq:fff-2sff} 
is definite. 
Then 
$f := \G \HH \G^*$, determined by \eqref{al:etaetabar}, 
\eqref{eq:mathcalG} and 
\eqref{eq:mtrixH}, is an HMC-1 immersion from $M$ to $\H^3$. 

Conversely, any HMC-1 surface (except a horosphere) 
has this parametrization in terms of $(G,h)$.     
\end{theorem}
\begin{remark}
\begin{enumerate}[(1)]
\item It has been already proved in \cite{KUY1} that 
 the solution to the differential equation \eqref{eq:g^-1dg} 
is described as \eqref{eq:mathcalG}.
 \item  Under the condition $h(z)=z$, the formula $f:= \G \HH \G^*$ 
with \eqref{eq:mathcalG} and \eqref{eq:mtrixH}, 
was already seen in \cite[Theorem 4]{GMM2}, where the condition $h(z)=z$ 
is caused by their assumption that $M$ is simply-connected and complete.   

As compared with it, Theorem \ref{thm:main1} is devoted to  
surfaces of non-trivial topology. The period condition is 
clarified, indeed, it is that $|\eta|^2$ is 
single-valued on $M$.   
\end{enumerate}
\end{remark}
The one-form $\theta$ can be calculated from the Schwarzian derivative as
follows:
\begin{lemma}
 \begin{equation*}
  \theta = - \frac{1}{2} \{ G; h \} dh 
\left( =
-\frac{1}{2} \left\{ \Bigl(\frac{G_{hh}}{G_h}\Bigr)_h%
-\frac{1}{2}%
\frac{(G_{hh})^2}{(G_h)^2}
\right\} dh
\right), 
 \end{equation*}
where $\{ G; h \}$ denotes the Schwarzian derivative of $G$
 with respect to $h$. 
\end{lemma}
\begin{proof}
Differentiating \eqref{eq:B^2} 
$B^2=-1/G_h$ with respect to $h$, we have 
\begin{equation*}
 dB/dh = \frac{1}{2} \frac{1}{B} \frac{G_{hh}}{(G_h)^2}. 
\end{equation*} 
Differentiating this again, we have
\begin{equation*}
 \frac{d}{dh} \left( \frac{dB}{dh} \right) 
=\frac{1}{2} \left\{ -\frac{1}{2}\frac{1}{B}%
\frac{(G_{hh})^2}{(G_h)^3}+\frac{1}{B}\frac{1}{G_h}%
\Bigl(\frac{G_{hh}}{G_h}\Bigr)_h\right\}. 
\end{equation*}
Therefore, it follows from \eqref{eq:theta} that  
\begin{equation*}
 \theta 
=\frac{1}{2} \left\{ -\frac{1}{2}\frac{1}{B^2}%
\frac{(G_{hh})^2}{(G_h)^3}+\frac{1}{B^2}\frac{1}{G_h}%
\Bigl(\frac{G_{hh}}{G_h}\Bigr)_h\right\} dh. 
\end{equation*}
Again, using $B^2=-1/G_h$, we obtain
\begin{equation*}
 \theta 
=\frac{1}{2} \left\{ \frac{1}{2}%
\frac{(G_{hh})^2}{(G_h)^2}-%
\Bigl(\frac{G_{hh}}{G_h}\Bigr)_h\right\} dh
= -\frac{1}{2} \{ G ; h \} dh.
\end{equation*}
\end{proof}

\section{Fronts with constant harmonic-mean curvature one}\label{sec:mrc-1front}
\subsection{Definition}
Let $M$ be a Riemann surface. 
Given a meromorphic function $G\colon M \to \C \cup \{ \infty \}$ 
and a pseudometric $|\eta|^2$ of constant curvature $-1$ on $M$,   
we can define a map $f = \G \HH \G^* \colon M \to \H^3$ 
 using \eqref{al:etaetabar}, \eqref{eq:mathcalG} and 
\eqref{eq:mtrixH}. Since $\G$ has poles $\{p_i\}$ in general, 
$f$ should be considered a map on $M \setminus \{p_i\}$. 
However, in such a case, we retake $M$ to be $M \setminus \{p_i\}$.  

We call $f$ an \emph{HMC-1 map} associated with $(G,|\eta|^2)$. 
By definition, the regular image of an HMC-1 map forms 
an immersed surface with constant harmonic-mean curvature one 
whose unit normal vector field is 
$\nu = \G 
\left[ \begin{smallmatrix}
 1 & -{\bar h} \\
 -{h} & {-1+|h|^2}
\end{smallmatrix} \right]
\G^*$.  
Though $f$ may fail to be an immersion, 
the unit normal $\nu$ is defined across the singularities. Hence, 
the following definition does make sense.

An HMC-1 map $f$ is called an \emph{HMC-1 front} if 
$(f, \nu) \colon M \to T_1 \H^3 (\cong T_1^* \H^3)$ is an immersion, 
where $T_1 \H^3$ ($T_1^* \H^3$) denotes the unit (co)tangent bundle
over $\H^3$. (The term {\it front} comes from {\it wave fronts} in 
the theory of singularities.)
It is obvious from the definition that the formulas for 
HMC-1 surfaces in the previous section 
can be applied for HMC-1 fronts.    
\begin{proposition}\label{prop:front}
For an HMC-1 map $f \colon M \to \H^3$, the following three conditions 
are equivalent:
\begin{enumerate}[\rm (1)]
 \item $f$ is an HMC-1 front. 
 \item The $(1,1)$-part $\fff^{1,1} (=|\eta|^2 + |\pi|^2)$ of the 
first fundamental form $\fff = |\eta + \pi|^2$ is 
a Riemannian metric on $M$. 
 \item $\G \colon {\tilde M} \to \SL(2, \C)$ is non-singular. 
\end{enumerate} 
\end{proposition}
\begin{proof}
We can put (1)--(3) in different words as follows: 
\begin{enumerate}[(1)]
 \item $|df|^2$ and $|d \nu|^2$ never vanish simultaneously, that is, 
$\bigl| {2 dh}/(1-|h|^2) +(1-|h|^2) \bar \theta \bigr|^2$ 
and $(1-|h|^2)^2 \left| \theta \right|^2$ never vanish simultaneously. 
 \item $\fff^{1,1} = {4 |dh|^2}/{(1-|h|^2)^2} + 
(1-|h|^2)^2 \left| \theta \right|^2$ never vanishes.
 \item Either $\theta \ne 0$,  or both $\theta =0$ and $dh \ne 0$.  
\end{enumerate} 
Then, it is not difficult to see the equivalency.
\end{proof}
\begin{remark}
 $T_1 \H^3$ ($\cong T_1^* \H^3$) has a canonical Riemannian metric, 
which is called the {\it Sasakian metric\/}. 
We denote by $\fff^S$ the pull-back 
of the Sasakian metric via the map $(f, \nu)$. $\fff^S$ is a 
Riemannian metric on $M$ for a front $f \colon M \to \H^3$. Indeed, 
$\fff^S = |df|^2+|d \nu|^2=|\eta + \pi|^2+|\pi|^2$. 
$\fff^S$ is not conformally equivalent to $\fff^{1,1}$ in general.      
\end{remark}

It is clear from Proposition \ref{prop:front} that 
a singularity of an HMC-1 front is a point where 
\begin{equation*}
 (\eta + \pi) \wedge (\bar \eta + \bar \pi)=0 \iff 
\eta \wedge \bar \eta + \pi \wedge \bar \pi =0 \iff
|\eta|^2 = |\pi|^2. 
\end{equation*}
\begin{proposition}
There are no compact HMC-1 fronts. 
\end{proposition}
\begin{proof}
 Suppose, by way of contradiction, that there exists a 
compact HMC-1 front $f \colon M \to \H^3$. 

It follows from Lemma \ref{lem:barpdpd} that 
\begin{equation*}
 (e_0+e_3)_{z \zb} = \frac{|\eta/dz|^2}{2} e_0. 
\end{equation*}
Taking the trace of both side, we have
 \begin{equation}\label{eq:subharmonic}
(\trace (e_0+e_3))_{z \zb} = \frac{|\eta/dz|^2}{2} \trace e_0 
= |\eta/dz|^2 x_0 \ge 0.  
\end{equation}
Hence, $\trace (e_0+e_3)$ is a subharmonic function on $M$. 
It must be constant, since $M$ is compact.  
Again by \eqref{eq:subharmonic}, we have $\eta = 0$, 
a contradiction.  
\end{proof}
\subsection{Weak completeness}
 We say that an HMC-1 front $f \colon M \to \H^3$ is 
\emph{weakly complete} if $\fff^{1,1}$, 
the $(1,1)$-part the first fundamental form $\fff$,  
is a complete Riemannian metric on $M$ (cf. \cite{KRSUY}). 
\begin{proposition}
 For an HMC-1 front, weak completeness is equivalent to 
the completeness of $\fff^S$, the induced metric of the Sasakian metric.   
\end{proposition}
\begin{proof}
Let $\gamma \colon [0, \infty) \to M$ be an arbitrary divergent path. 
Recall that $\fff^{S}= |\eta + \pi |^2 + |\pi|^2$ and 
$\fff^{1,1}=|\eta |^2 + |\pi|^2$.   
We wish to prove that, if $\gamma$ has infinite length with respect to one 
of the two metrics, then it also has infinite length 
with respect to the other metric.  
When $\int_{\gamma} |\pi| = \infty$, it is trivial that 
$\gamma$ has infinite length with respect to both metrics 
$\fff^S$ and $\fff^{1,1}$. 
Hence, we have only to give a proof under the 
assumption $\int_{\gamma} |\pi| < \infty$. 
\begin{enumerate}[(a)]
 \item Suppose that $\fff^{1,1}$ is complete. Clearly, 
$\int_{\gamma} |\eta| = \infty$. If we divide the interval $[0, \infty)$ 
so that
 \begin{equation*}
 [0, \infty)=J_+ \cup J_-, \text{ where } 
J_+ = \{ |\eta| \ge |\pi| \}, \ J_- = \{ |\eta| < |\pi| \}, 
\end{equation*} 
then $\int_{J_+} |\eta| = \infty$, because 
$\int_{J_-} |\eta| < \int_{J_-} |\pi| < \infty$. Thus, 
\begin{align*}
  &  \int_{\gamma} \sqrt{|\eta + \pi|^2+|\pi|^2} 
\ge \int_{\gamma} |\eta+\pi|
\ge \int_{\gamma} \left| |\eta|-|\pi| \right| \\
\ge & \int_{J_+} |\eta|-|\pi| 
= \int_{J_+} |\eta|- \int_{J_+}|\pi|
=\infty - \text{(finite value)} = \infty. 
\end{align*}
Therefore, $\fff^S$ is complete. 
\item Conversely, we suppose $\fff^S$ is complete. Clearly, 
$\int_{\gamma} |\eta+\pi| = \infty$. If we divide the interval $[0, \infty)$ 
so that  
\begin{equation*}
 [0, \infty)=J'_+ \cup J'_-, \text{ where } 
J'_+ = \{ |\eta+\pi| \ge |\pi| \}, \ J'_- = \{ |\eta+\pi| < |\pi| \}, 
\end{equation*}
then $\int_{J'_+} |\eta+\pi| = \infty$, because 
$\int_{J'_-} |\eta+\pi| < \int_{J'_-} |\pi| = \infty$. Thus 
\begin{align*}
  &  \int_{\gamma} \sqrt{|\eta|^2+|\pi|^2} 
\ge \int_{\gamma} |\eta|
\ge \int_{\gamma} \left| |\eta+\pi|-|-\pi| \right| \\
\ge & \int_{J'_+} |\eta+\pi|-|\pi| 
\ge \int_{J'_+} |\eta+\pi| - \int_{J'_+} |\pi|= \infty. 
\end{align*}
Therefore $\fff^{1,1}$ is complete. 
\end{enumerate}
\end{proof}
 
Note that an HMC-1 front is weakly complete if 
it is complete (in the usual sense),  
because $\fff^S = |df|^2 + |d \nu|^2$ is complete if 
$\fff = |df|^2$ is complete. 
\subsection{HMC-1 fronts of finite topology}
 There are two kind of ends for (weakly) complete HMC-1 fronts 
of finite topology. One is conformally equivalent to the punctured disk 
$\varDelta^* = \setdef{z}{0 < |z| <1}$, and the other is conformally 
equivalent to the annulus $A_r=\setdef{z}{r < |z| <1}$. We shall 
call the former a {\it puncture-type end\/}, the latter an 
{\it annular end\/}. 
For a puncture-type end $\varDelta^*$, we also call a point $z=0$ 
an end. 
For an annular end $A_r$, we also call the boundary $|z|=r$ 
an end. 
\begin{theorem}\label{prop:p-end}
Let $f \colon M \to \H^3$ be 
a weakly complete HMC-1 front of finite topology. Then  
\begin{enumerate}[\rm (1)]
 \item  the set of singularities never accumulate to 
a puncture-type end, and 
 \item the Gaussian curvature $K(z)$ converges to $0$ as 
$z$ tends to a puncture-type end. Hence, the mean curvature $H(z)$
converges to $1$.   
\end{enumerate}
\end{theorem}
\begin{proof}
 Let $ \varDelta^* \subset M$ be a puncture-type end. 
If we assumed $\lim_{z \to 0} h(z) \in \pd \D$, then some 
portion of a neighborhood of $0$ is not contained in the image of 
$h$. This is a contradiction. Therefore $h(0) \in \D$. Hence, 
$\lim_{z \to 0} |\eta|^2
=\lim_{z \to 0}{4|dh|^2}/{(1-|h|^2)^2} < \infty$. 

On the other hand,  
$\fff^{1,1}= |\eta|^2 + |\pi|^2=
{4|dh|^2}/{(1-|h|^2)^2}+(1-|h|^2)^2|\theta|^2$ is 
complete at $0$. Thus 
$\lim_{z \to 0} |\pi|^2= \infty$. 

Therefore, $|\pi|^2 \ne |\eta|^2$ near $0$,  
and  
$\lim_{z \to 0} K = 
\lim_{z \to 0} {|\eta|^2}/{ (|\pi|^2 - |\eta|^2)} =0$. 
\end{proof}
Note that, compared with Proposition \ref{prop:p-end} (1), 
it can occur that the set of singularities accumulate 
toward an annular end (see Example \ref{eg:z+1/z,z} below).
\begin{proposition}\label{prop:a-end}
Let $f \colon M \to \H^3$ be 
a weakly complete HMC-1 front of finite topology, and 
$V ( \cong A_r)$ an annular end. 
 Then the Gaussian curvature $K(z)$ converges to $-1$ as 
$z$ tends to a point $z_0 \in \pd A_r$ with $|z_0|=r$, 
unless $\lim_{z \to z_0} \tff = \infty$.  
Hence, the mean curvature $H(z)$ converges to $0$.   
\end{proposition}
\begin{proof}
Let $A_r \subset M$ be an annular end, and 
 take an arbitrary point $z_0 \in \pd A_r$ with $|z_0|=r$. 
Then $\lim_{z \to z_0} h \in \pd \D$. 
(If we assumed $\lim_{z \to z_0} h \in \D$, 
then the image of some neighborhood of $z_0$ 
is also contained in $\D$. At every point $w$ in 
the neighborhood, $\lim_{z \to w}|\theta|$ must be infinity, 
since all divergent paths have infinite length. However, this is impossible 
because $\theta$ is a holomorphic one-form.) 
Hence, $\lim_{z \to z_0} |\eta|^2=\infty$. 
Therefore, 
$\lim_{z \to z_0} K = 
\lim_{z \to z_0} {|\eta|^2}/{(|\pi|^2 - |\eta|^2)} =-1$, 
unless $\lim_{z \to z_0} \tff =\lim_{z \to z_0} |\pi|^2 = \infty$.  
\end{proof}
From Proposition \ref{prop:p-end} and Proposition \ref{prop:a-end}, 
we may use the adjective \emph{horospherical} for a puncture-type end,   
and \emph{hemispherical} for an annular end. 

\subsection{Examples}
We show some examples of (weakly) complete HMC-1 fronts 
of finite topology. 

In this section, we denote by $ds^2_H$ the Poincar\'e metric 
on the unit disk, and in the figures, the hyperbolic three-space 
$\H^3$ is realized by the Poincar\'e ball model. 

\begin{example}
 For a positive number $\alpha$, consider 
\begin{equation*}
 G(z)=z, \ |\eta|^2=\frac{4 |\alpha|^2 |z|^{2\alpha-2}}%
{(1-|z|^{2 \alpha})^2}|dz|^2
\left(= h^{*} ds^2_H \text{ where } h(z)=z^{\alpha} \right)
\end{equation*}
on 
\begin{equation*}
 M = \begin{cases}
      \varDelta=\setdef{z}{|z|<1} & \text{if $\alpha=1$, } \\
      \varDelta \setminus{ \{ 0 \} } & \text{otherwise. }
     \end{cases}
\end{equation*}
Then the HMC-1 front $f \colon M \to \H^3$ 
associated with $(G,|\eta|^2)$ satisfies 
\begin{align*}
& \theta = \frac{1-\alpha^2}{4 \alpha} z^{-\alpha-1}dz, \ 
dh= \alpha z^{\alpha-1} dz, \ 
Q(=\eta \bar \pi) = \frac{1-\alpha^2}{2z^2} dz^2, \\
& |\pi|^2=\frac{|\alpha^2 -1|^2}{16 |\alpha|^2}
\frac{(1-|z|^{2 \alpha})^2}{|z|^{2\alpha+2}}|dz|^2. 
\end{align*}
If $\alpha=1$, then $f$ is complete and totally geodesic.    
If $\alpha \ne 1$, then $f$ is weakly complete and has 
a horospherical end at $z=0$ and a hemispherical end at 
$|z|=1$. 
Its singular locus is the circle  
\begin{equation*}
 |z| = \left(- \sqrt{\frac{2 \alpha^2}{\alpha^2-1}}
+ \sqrt{\frac{2 \alpha^2}{\alpha^2-1} +1} \right)^{1/\alpha}. 
\end{equation*}
$K \ge 0$ inside this circle, and $K \le -1$ outside the circle. 
\end{example}
\begin{figure}[ht]
\begin{center}
 \begin{tabular}{cc}
\includegraphics[width=4.5cm,clip,keepaspectratio]{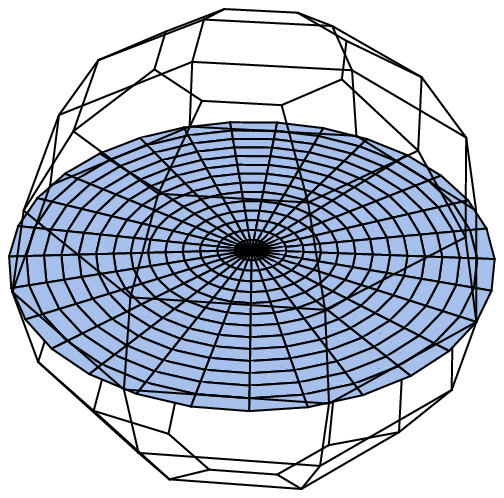}&
\includegraphics[width=4.5cm,clip,keepaspectratio]{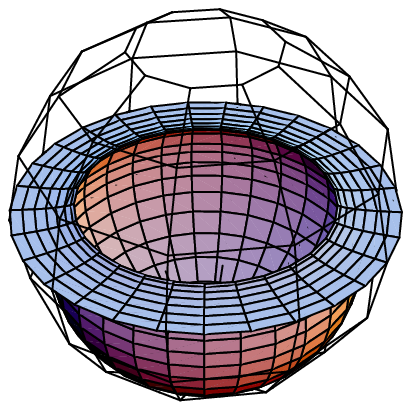} \\
\small $(G,h)=(z,z)$ &
 \small $(G,h)=(z,z^2)$ \\
\includegraphics[width=4.5cm,clip,keepaspectratio]{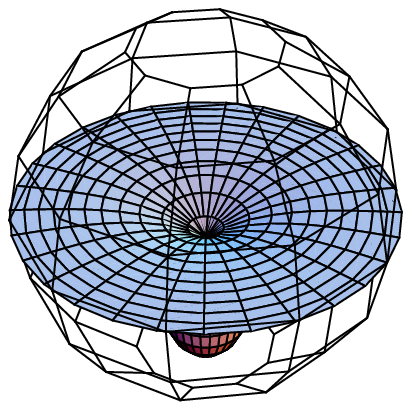}&
\includegraphics[width=4.5cm,clip,keepaspectratio]{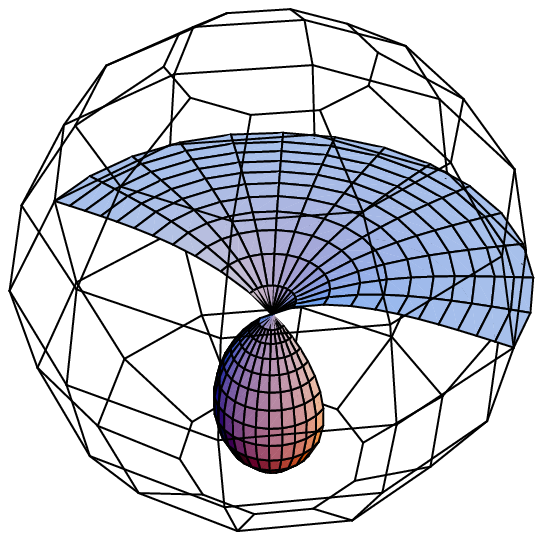}\\
\small $(G,h)=(z,\sqrt{z})$ &
 \small $(G,h)=(z, \sqrt{z})$ (half cut)
\end{tabular}
\end{center}
\caption{}\label{fig:hemi}
\end{figure}
\begin{example}
For a non-zero complex number $k$, consider 
 $G(z)=\exp (kz) , \ |\eta|^2 = ds^2_H$ on $\varDelta$, 
then the associated HMC-1 front satisfies 
\begin{align*}
  \theta =\frac{k^2}{4}dz, \ dh=dz, \ 
 Q = \frac{k^2}{2} dz^2, \ 
|\pi|^2=\frac{|k|^4}{16}(1-|z|^2)^2 |dz|^2.  
\end{align*}
Its singular locus is $|z|^2=1-{2 \sqrt{2}}/{|k|}$. 
In particular, it has no singularities if $|k| < 2 \sqrt{2}$. 
\end{example}
\begin{figure}[ht]
\begin{center}
 \begin{tabular}{ccc}
\includegraphics[width=4.5cm,clip,keepaspectratio]{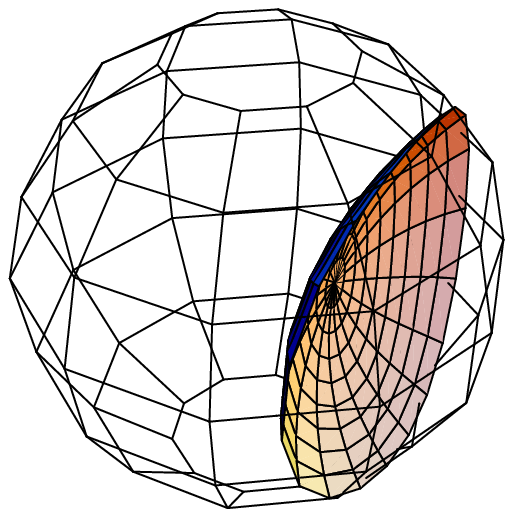}&
\includegraphics[width=4.5cm,clip,keepaspectratio]{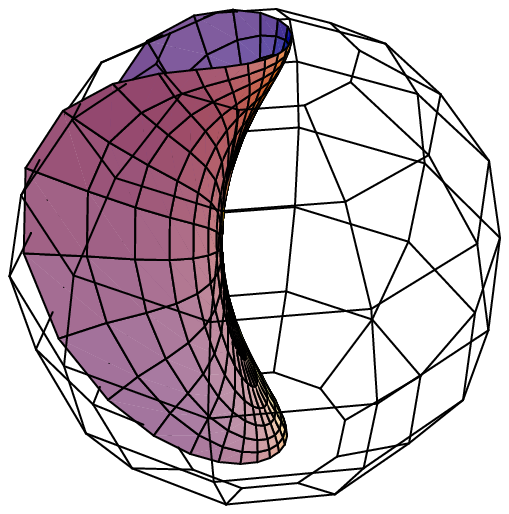}&
\includegraphics[width=4.5cm,clip,keepaspectratio]{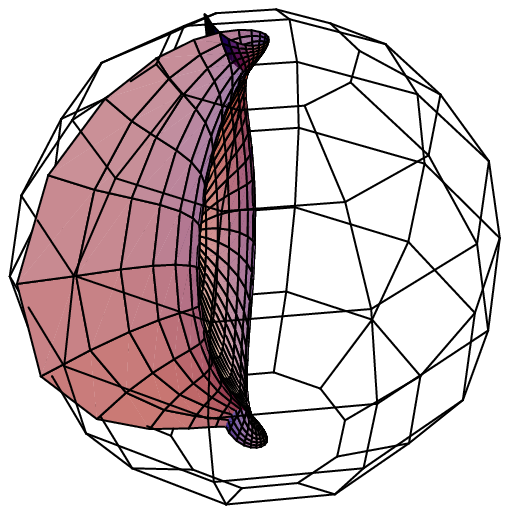}\\
\small $(G,h)=(\exp z,z)$ &
 \small $(G,h)=(\exp 2\sqrt{2} z,z)$ &
\small $(G,h)=(\exp 4z,z)$ 
\end{tabular}
\end{center}
\caption{}
\end{figure}
%
%
\begin{example}\label{eg:z+1/z,z}
Consider $G(z)=z+ {1}/{z}, \ |\eta|^2 = ds^2_H$ on $\varDelta$. 
Then the associated HMC-1 front satisfies 
\begin{align*}
  \theta =\frac{3}{(z^2-1)^2}dz, \ dh=dz, \ 
Q = \frac{6}{(z^2-1)^2} dz^2, \ 
|\pi|^2= \frac{9 (1-|z|^2)^2}{|z^2-1|^4} |dz|^2.  
\end{align*}
Its singular locus is 
\begin{equation*}
C : 2 |z^2-1|^2 = 3 (1-|z|^2)^2. 
\end{equation*} 
$f$ has an annular end at $|z|=1$, 
and  the singular locus $C$ accumulates at $z= \pm 1$.  
\end{example}
\begin{figure}[ht]
\begin{center}
 \begin{tabular}{cc}
\includegraphics[width=4.5cm,clip,keepaspectratio]{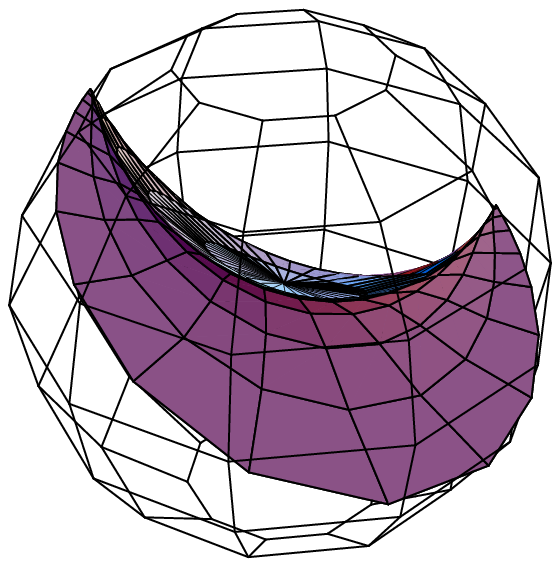}&
\includegraphics[bb = 91 -83 322 335, width=3.0cm,clip,keepaspectratio]{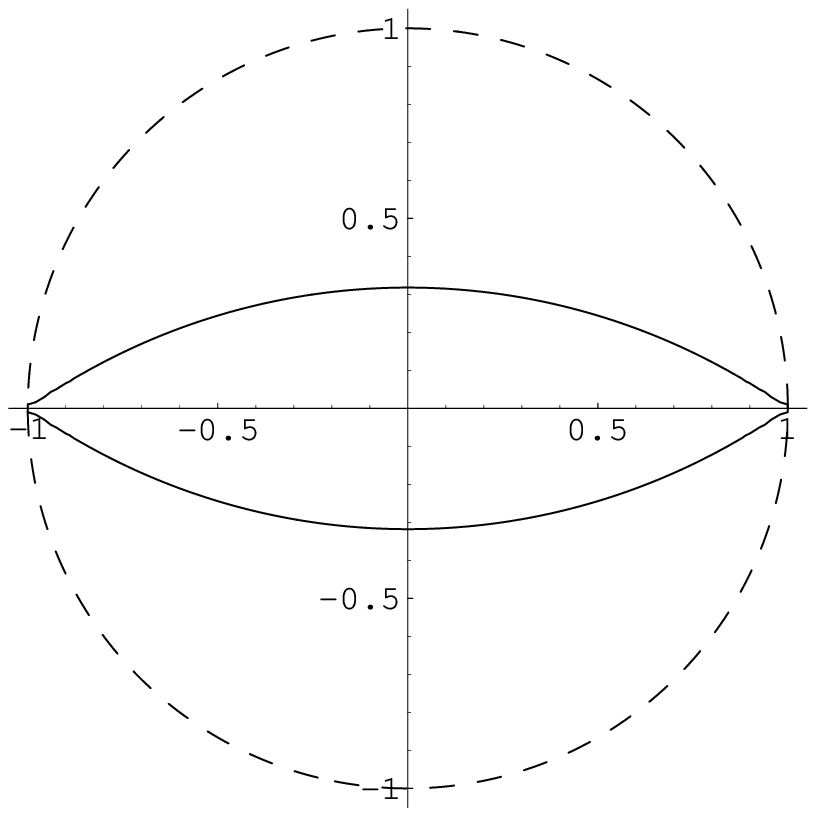}\\
\small $(G,h)=(z+1/z,z)$ &
 \small singular locus in $\varDelta$
\end{tabular}
\end{center}
\caption{}
\end{figure}



\end{document}